\documentclass[letterpaper,twocolumn,conference]{IEEEtran}
\usepackage[T1]{fontenc}
\usepackage[latin9]{inputenc}
\usepackage{units}
\usepackage{amsmath}
\usepackage{amsthm}
\usepackage{amssymb}

\makeatletter

\pdfpageheight\paperheight
\pdfpagewidth\paperwidth

\theoremstyle{plain}
\newtheorem{thm}{\protect\theoremname}
\theoremstyle{remark}
\newtheorem{rem}[thm]{\protect\remarkname}
\theoremstyle{plain}
\newtheorem{lem}[thm]{\protect\lemmaname}
\theoremstyle{plain}
\newtheorem{cor}[thm]{\protect\corollaryname}

\usepackage[caption=false,font=footnotesize]{subfig}

\makeatother

\providecommand{\corollaryname}{Corollary}
\providecommand{\lemmaname}{Lemma}
\providecommand{\remarkname}{Remark}
\providecommand{\theoremname}{Theorem}

\begin{document}

\title{Statistical Inference and Exact Saddle Point Approximations}

\author{\IEEEauthorblockN{Peter~Harremo{\"e}s}\IEEEauthorblockA{GSK Department\\
Copenhagen Business College\\
Copenhagen, Denmark\\
Email: harremoes@ieee.org}}
\maketitle
\begin{abstract}
Statistical inference may follow a frequentist approach or it may
follow a Bayesian approach or it may use the minimum description length
principle (MDL). Our goal is to identify situations in which these
different approaches to statistical inference coincide. It is proved
that for exponential families MDL and Bayesian inference coincide
if and only if the renormalized saddle point approximation for the
conjugated exponential family is exact. For 1-dimensional exponential
families the only families with exact renormalized saddle point approximations
are the Gaussian location family, the Gamma family and the inverse
Gaussian family. They are conjugated families of the Gaussian location
family, the Gamma family and the Poisson-exponential family. The first
two families are self-conjugated implying that only for the two first
families the Bayesian approach is consistent with the frequentist
approach. In higher dimensions there are more examples.
\end{abstract}

\section{Introduction}

In this paper we are interested in predition of the future given the
past. We assume that a sequence $x_{1}^{m}=x_{1},x_{2},\dots,x_{m}$
has been observed and the goal is to predict the next symbols $x_{m+1}^{n}=x_{m+1},x_{m+2},\dots,x_{n}$
in the sense that we will assign a probability or a probability density
to this sequence. The prediction is compared with iid models given
by a parametrized family $\left(P_{\theta}\right)_{\theta\in\Theta}$
of probability distributions that assign probability $\prod_{i=1}^{n}P_{\theta}\left(x_{i}\right)$
(or the corresponding density) to the sequence $x_{1}^{n}.$ One may
think of the elements of the family $\left(P_{\theta}\right)_{\theta\in\Theta}$
as the models that some experts can choose among. For the techniques
used in this paper the restriction to iid models is crusial, but some
of the results may generalize to non-iid models.

All measures will be described by their density with respect to a
dominating measure $\lambda.$ Data are assumed to lie in $\mathcal{X}\subseteq\mathbb{R}^{d}$
and vectors will be marked with bold face. Assume that $\left(P_{\boldsymbol{\theta}}\right)_{\boldsymbol{\theta}\in\Theta}$
is a natural exponential family with 
\begin{align*}
\frac{\mathrm{d}P_{\boldsymbol{\theta}}}{\mathrm{d}\lambda}\left(\boldsymbol{x}\right) & =\frac{\exp\left(\boldsymbol{\theta}\cdot\boldsymbol{x}\right)}{Z\left(\boldsymbol{\theta}\right)}=\exp\left(\boldsymbol{\theta}\cdot\boldsymbol{x}-A\left(\boldsymbol{\theta}\right)\right).
\end{align*}
Here $Z\left(\boldsymbol{\theta}\right)=\int\exp\left(\boldsymbol{\theta}\cdot\boldsymbol{x}\right)\,\mathrm{d}\lambda\boldsymbol{x}$
is the \emph{moment generating function} and $A\left(\boldsymbol{\theta}\right)=\ln\left(Z\left(\boldsymbol{\theta}\right)\right)$
is the \emph{cumulant generating function}. If the parameter has value
$\boldsymbol{\theta}$ then the mean value is $\boldsymbol{\mu}_{\boldsymbol{\theta}}=\nabla A\left(\boldsymbol{\theta}\right).$
The density $\frac{\mathrm{d}P_{\boldsymbol{\theta}}}{\mathrm{d}\lambda}$
will be denoted $p_{\boldsymbol{\theta}}$, but sometimes we will
also use $p_{\boldsymbol{\theta}}$ for iid sequences. 

One approach is the frequentist approach where the sequence $\boldsymbol{x}_{1}^{n}$
is generated by the distribution $P_{\boldsymbol{\theta}}$ for some
true but unknown value of $\boldsymbol{\theta}.$ The sequence $\boldsymbol{x}_{1}^{m}$
is used to make inference about the value of $\boldsymbol{\theta}$
in terms of a confidence region. In a Bayesian approach one has a
prior distribution $\pi$ on the true parameter $\boldsymbol{\theta}$
and the sequence $\boldsymbol{x}_{1}^{m}$ is used to calculate a
posterior distribution of $\boldsymbol{\theta}$ as 
\[
\frac{p_{\boldsymbol{\theta}}\left(\boldsymbol{x}^{m}\right)\pi\left(\boldsymbol{\theta}\right)}{\int_{\Theta}p_{\boldsymbol{\theta}}\left(\boldsymbol{x}^{m}\right)\pi\left(\boldsymbol{\theta}\right)\,\mathrm{d}\boldsymbol{\theta}}.
\]
Then the posterior distribution of $\boldsymbol{x}_{m+1}^{n}$ is
given by 
\begin{align}
p_{\pi}\left(\boldsymbol{x}_{m+1}^{n}\mid\boldsymbol{x}^{m}\right) & =\int_{\Theta}p_{\boldsymbol{\theta}}\left(\boldsymbol{x}_{m+1}^{n}\right)\,\mathrm{d}\pi\left(\boldsymbol{\theta}\mid\boldsymbol{x}^{m}\right)\label{eq:posteriorpredic}\\
 & =\int_{\Theta}p_{\boldsymbol{\theta}}\left(\boldsymbol{x}_{m+1}^{n}\right)\frac{p_{\boldsymbol{\theta}}\left(\boldsymbol{x}^{m}\right)\pi\left(\boldsymbol{\theta}\right)}{\int_{\Theta}p_{\boldsymbol{\theta}}\left(\boldsymbol{x}^{m}\right)\pi\left(\boldsymbol{\theta}\right)\,\mathrm{d}\boldsymbol{\theta}}\,\mathrm{d}\boldsymbol{\theta}\nonumber 
\end{align}
One of the main problems is Bayesian statistics is the question of
how to determine the prior distribution $\pi.$ 

The moment generating function $Z$ is related to the Laplace transform
of the measure $\lambda$, so any of the functions $Z$ and $A$ can
be used to reconstruct $\lambda$. The \emph{Hesse matrix} of $A$
with respect to $\boldsymbol{\theta}$ equals the \emph{co-variance
matrix} $Cov\left(\boldsymbol{\mu}_{\boldsymbol{\theta}}\right)$.
The Fisher information matrix with respect to the natural parameter
is $Cov\left(\boldsymbol{\mu}_{\boldsymbol{\theta}}\right)$ so that
\emph{Jeffreys' prior} is proportional to $\left|Cov\left(\boldsymbol{\mu}_{\boldsymbol{\theta}}\right)\right|^{\nicefrac{1}{2}}.$
Therefore \emph{Jeffreys' posterior} distribution of the parameter
$\boldsymbol{\theta}$ after observing a sequence of length $m$ with
average $\bar{\boldsymbol{x}}$ is proportional to
\[
\exp\left(m\cdot\left(\boldsymbol{\theta}\cdot\bar{\boldsymbol{x}}-A\left(\boldsymbol{\theta}\right)\right)\right)\cdot\left|Cov\left(\boldsymbol{\mu}_{\boldsymbol{\theta}}\right)\right|^{\nicefrac{1}{2}}.
\]
One motivation for using Jeffreys' prior is that it is considered
as an uninformative prior. Another motivation is that if one restricts
to a bounded subset whose closure is in the interior of the full parameter
space, then the use of Jeffrey's prior is asymptotically optimal in
a MDL sense \cite{Grunwald2007}.

A co-variance matrix is positive semi-definite so the cumulant generating
function is convex. The\emph{ convex conjugate} of the cumulant generating
function $A$ is $A^{*}\left(\boldsymbol{x}\right)=\sup_{\boldsymbol{\theta}}\left\{ \boldsymbol{\theta}\cdot\boldsymbol{x}-A\left(\boldsymbol{\theta}\right)\right\} .$
The conjugate parameter $\boldsymbol{x}^{*}$ equals the value of
$\boldsymbol{\theta}$ such that $P_{\boldsymbol{\theta}}$ has mean
value $\boldsymbol{x}$, i.e. $\boldsymbol{x}^{*}$ is the solution
to the equation $\nabla A\left(\boldsymbol{\theta}\right)=\boldsymbol{x}$.
Usually the conjugate parameter $\boldsymbol{x}^{*}$ is denoted $\hat{\boldsymbol{\theta}}\left(\boldsymbol{x}\right)$
and is called the maximum likelihood estimate of $\boldsymbol{\theta}.$
We can define the \emph{conjugated exponential family} (if it exists)
as the exponential family with sufficient statistic $\boldsymbol{\theta}$
and with cumulant generating function $A^{*}\left(\boldsymbol{x}\right).$ 
\begin{rem}
For an exponential family the conjugated exponential family gives
a set of ``conjugated priors'' as this concept is defined in the
literature on Bayesian statistics (see \cite{Raiffa1961} and \cite[Sec. 12.2.6]{Liu2014}),
but a set of ``conjugated priors'' need not coincide with the conjugated
exponential family as it is defined in this paper.
\end{rem}
The Bregman divergence generated by the convex function $A$ is defined
by
\[
D_{A}\left(\boldsymbol{\theta}_{2},\boldsymbol{\theta}_{1}\right)=A\left(\boldsymbol{\theta}_{2}\right)-\left(A\left(\boldsymbol{\theta}_{1}\right)+\left(\boldsymbol{\theta}_{2}-\boldsymbol{\theta}_{1}\right)\cdot\nabla A\left(\boldsymbol{\theta}_{1}\right)\right)
\]
 Using convex conjugation the divergence can also be written as
\[
D_{A}\left(\boldsymbol{\theta}_{2},\boldsymbol{\theta}_{1}\right)=D_{A^{*}}\left(\boldsymbol{\mu}_{1},\boldsymbol{\mu}_{2}\right).
\]
The information divergence can be calculated as
\begin{multline*}
D\left(\left.P_{\boldsymbol{\theta}_{1}}\right\Vert P_{\boldsymbol{\theta}_{2}}\right)=E_{\boldsymbol{\theta}_{1}}\left[\ln\left(\frac{\mathrm{d}P_{\boldsymbol{\theta}_{1}}}{\mathrm{d}P_{\boldsymbol{\theta}_{2}}}\right)\right]\\
=E_{\boldsymbol{\theta}_{1}}\left[\left(\boldsymbol{\theta}_{1}\cdot\boldsymbol{X}-A\left(\boldsymbol{\theta}_{1}\right)\right)-\left(\boldsymbol{\theta}_{2}\cdot\boldsymbol{X}-A\left(\boldsymbol{\theta}_{2}\right)\right)\right]\\
=D_{A}\left(\boldsymbol{\theta}_{2},\boldsymbol{\theta}_{1}\right).
\end{multline*}

The conjugated exponential family gives posterior distributions on
the parameter $\boldsymbol{\theta},$ such that the maximum likelihood
estimate $\hat{\boldsymbol{\theta}}\left(\boldsymbol{x}\right)$ is
unbiased in the sense that it equals the mean value of $\boldsymbol{\theta}$
with respect to the posterior distribution of $\boldsymbol{\theta}$
given $\boldsymbol{x}$. Therefore the use of the conjugated exponential
family implies that the maximum likelihood estimator equals the Bayes
estimator with respect to the loss function $D_{A}$ or any other
Bregman divergence. 

The likelihood function can be written as
\begin{multline*}
p_{\boldsymbol{\theta}}\left(\boldsymbol{x}\right)=\exp\left(\boldsymbol{\theta}\cdot\boldsymbol{x}-A\left(\boldsymbol{\theta}\right)\right)=\\
\exp\left(-A\left(\boldsymbol{\theta}\right)+\left(A\left(\hat{\boldsymbol{\theta}}\left(\boldsymbol{x}\right)\right)+\left(\boldsymbol{\theta}-\hat{\boldsymbol{\theta}}\left(\boldsymbol{x}\right)\right)\cdot\nabla A\left(\hat{\boldsymbol{\theta}}\left(\boldsymbol{x}\right)\right)\right)\right)\\
\cdot p_{\hat{\boldsymbol{\theta}}\left(\boldsymbol{x}\right)}\left(\boldsymbol{x}\right)\\
=\exp\left(-D_{A}\left(\boldsymbol{\theta},\hat{\boldsymbol{\theta}}\left(\boldsymbol{x}\right)\right)\right)\cdot p_{\hat{\boldsymbol{\theta}}\left(\boldsymbol{x}\right)}\left(\boldsymbol{x}\right).
\end{multline*}
As a consequence we have the following robustness property \cite[Section 19.3, Eq. 19.12]{Grunwald2007}
of the exponential family
\begin{equation}
\frac{\mathrm{d}P_{\boldsymbol{\theta}}}{\mathrm{d}P_{\hat{\boldsymbol{\theta}}\left(\boldsymbol{x}\right)}}\left(\boldsymbol{x}\right)=\exp\left(-D_{A}\left(\boldsymbol{\theta},\hat{\boldsymbol{\theta}}\left(\boldsymbol{x}\right)\right)\right).\label{eq:robust}
\end{equation}
The likelihood function after observing the sequence $\boldsymbol{x}^{m}$
is 
\begin{multline*}
\prod_{i=1}^{m}p_{\boldsymbol{\theta}}\left(\boldsymbol{x}_{i}\right)=\prod_{i=1}^{m}\exp\left(\boldsymbol{\theta}\cdot\boldsymbol{x}_{i}-A\left(\boldsymbol{\theta}\right)\right)\\
=\exp\left(\boldsymbol{\theta}\cdot\sum_{i=1}^{m}\boldsymbol{x}_{i}-m\cdot A\left(\boldsymbol{\theta}\right)\right)\\
=\exp\left(m\cdot\left(\boldsymbol{\theta}\cdot\bar{\boldsymbol{x}}-A\left(\boldsymbol{\theta}\right)\right)\right)\\
=\exp\left(-m\cdot D_{A}\left(\boldsymbol{\theta},\hat{\boldsymbol{\theta}}\left(\bar{\boldsymbol{x}}\right)\right)\right)\\
\cdot\exp\left(m\left(\hat{\boldsymbol{\theta}}\left(\bar{\boldsymbol{x}}\right)\cdot\bar{\boldsymbol{x}}-A\left(\hat{\boldsymbol{\theta}}\left(\bar{\boldsymbol{x}}\right)\right)\right)\right).
\end{multline*}

In the minimum description length (MDL) approach to statistical inference
there is no assumption about a true value of $\boldsymbol{\theta}$,
and the quality of a prediction is compared with the maximum likelihood
estimate of $\boldsymbol{\theta}$ in terms of a difference in code
length. For a data sequence $\boldsymbol{x}^{n}$ the \emph{regret}
of predicting $p\left(\boldsymbol{x}_{m+1}^{n}\mid\boldsymbol{x}^{m}\right)$
is
\[
-\ln\left(p\left(\boldsymbol{x}_{m+1}^{n}\mid\boldsymbol{x}^{m}\right)\right)-\left(-\ln\left(p_{\hat{\boldsymbol{\theta}}\left(\boldsymbol{x}^{n}\right)}\left(\boldsymbol{x}^{n}\right)\right)\right).
\]
Here the predictor $p\left(\cdot\mid\boldsymbol{x}^{m}\right)$ is
used to code the future $\boldsymbol{x}_{m+1}^{n}$ while the expert
is coding the whole sequence $\boldsymbol{x}^{n}$, but the expert
is allowed to choose the model $\boldsymbol{\theta}=\hat{\boldsymbol{\theta}}\left(\boldsymbol{x}^{n}\right)$
that gives the best fit to data. We take the maximum over all possible
data sequences and the predictor that minimizes the maximal regret
is called the \emph{conditional normalized maximum likelihood} predictor
(CNML) \cite{Rissanen2007} and is given by
\begin{multline}
p_{cnml}^{n}\left(\boldsymbol{x}_{m+1}^{n}\mid\boldsymbol{x}^{m}\right)\\
=\frac{p_{\hat{\boldsymbol{\theta}}\left(\boldsymbol{x}^{n}\right)}\left(\boldsymbol{x}^{n}\right)}{\int_{\mathcal{X}^{n-m}}p_{\hat{\boldsymbol{\theta}}\left(\boldsymbol{x}^{m}\boldsymbol{y}^{n-m}\right)}\left(\boldsymbol{x}^{m}\boldsymbol{y}^{n-m}\right)\,\mathrm{d}\lambda^{n-m}\left(\boldsymbol{y}^{n-m}\right)}.\label{eq:CNML}
\end{multline}

\section{Main results}

The essence of the following lemma was already present in \cite[Lem. 3]{Bartlett2013}. 
\begin{lem}
\label{lem:key}Assume that $\left(P_{\boldsymbol{\theta}}\right)_{\boldsymbol{\theta}\in\Theta}$
is a natural exponential family. Assume that $m$ is a number such
that CNML and Bayesian prediction based on a prior $\pi$ give equal
prediction strategies for sequences $\boldsymbol{x}_{m+1}^{n}$ for
all $n>m.$ Then for any $n>m$ the integral
\[
\int_{\Theta}\frac{p_{\boldsymbol{\theta}}\left(\boldsymbol{x}^{n}\right)}{p_{\hat{\boldsymbol{\theta}}\left(\boldsymbol{x}^{n}\right)}\left(\boldsymbol{x}^{n}\right)}\pi\left(\boldsymbol{\theta}\right)\,\mathrm{d}\boldsymbol{\theta}
\]
is constant as a function of the data sequence $\boldsymbol{x}^{n}=\boldsymbol{x}_{1}\boldsymbol{x}_{2}\dots\boldsymbol{x}_{n}$~,
\end{lem}
\begin{rem}
Prediction with CNML and prediction based of Jeffreys prior can only
be equal if they are both defined. The values of $m$ for which these
prediction methods are defined, may in principle be different and
may depend on the data sequence \cite{Harremoes2013}. 
\end{rem}
\begin{IEEEproof}
For all $\boldsymbol{x^{n}}\in\mathcal{X}^{n}$ we must have
\[
p_{\pi}\left(\boldsymbol{x}_{m+1}^{n}\mid\boldsymbol{x}_{m}\right)=p_{cnml}^{n}\left(\boldsymbol{x}_{m+1}^{n}\mid\boldsymbol{x}_{m}\right).
\]
Using (\ref{eq:posteriorpredic}) and (\ref{eq:CNML}) we get
\begin{multline*}
\int_{\Theta}p_{\boldsymbol{\theta}}\left(\boldsymbol{x}_{m+1}^{n}\right)\frac{p_{\boldsymbol{\theta}}\left(\boldsymbol{x}^{m}\right)\pi\left(\boldsymbol{\theta}\right)}{\int_{\Theta}p_{\boldsymbol{\theta}}\left(\boldsymbol{x}^{m}\right)\pi\left(\boldsymbol{\theta}\right)\,\mathrm{d}\boldsymbol{\theta}}\,\mathrm{d}\boldsymbol{\theta}\\
=\frac{p_{\hat{\boldsymbol{\theta}}\left(\boldsymbol{x}^{n}\right)}\left(\boldsymbol{x}^{n}\right)}{\int_{\mathcal{X}^{n-m}}p_{\hat{\boldsymbol{\theta}}\left(\boldsymbol{x}^{m}\boldsymbol{y}^{n-m}\right)}\left(\boldsymbol{x}^{m}\boldsymbol{y}^{n-m}\right)\,\mathrm{d}\lambda^{n-m}\left(\boldsymbol{y}^{n-m}\right)}
\end{multline*}
and 
\begin{multline*}
\frac{\int_{\Theta}p_{\boldsymbol{\theta}}\left(\boldsymbol{x}^{n}\right)\pi\left(\boldsymbol{\theta}\right)\,\mathrm{d}\boldsymbol{\theta}}{p_{\hat{\boldsymbol{\theta}}\left(\boldsymbol{x}^{n}\right)}\left(\boldsymbol{x}^{n}\right)}\\
=\frac{\int_{\Theta}p_{\boldsymbol{\theta}}\left(\boldsymbol{x}^{m}\right)\pi\left(\boldsymbol{\theta}\right)\,\mathrm{d}\boldsymbol{\theta}}{\int_{\mathcal{X}^{n-m}}p_{\hat{\boldsymbol{\theta}}\left(\boldsymbol{x}^{m}\boldsymbol{y}^{n-m}\right)}\left(\boldsymbol{x}^{m}\boldsymbol{y}^{n-m}\right)\,\mathrm{d}\lambda^{n-m}\left(\boldsymbol{y}^{n-m}\right)}.
\end{multline*}
The quantity on the left side is a function of $\boldsymbol{x}^{n}$
while the quantity on the right side is a function of the sub-string
$\boldsymbol{x}^{m}.$ Since the model is invariant under permutations
of the elements in the string $\boldsymbol{x}^{n}$ both sides must
equal a constant. Finally we note that
\[
\frac{\int_{\Theta}p_{\boldsymbol{\theta}}\left(\boldsymbol{x}^{n}\right)\pi\left(\boldsymbol{\theta}\right)\,\mathrm{d}\boldsymbol{\theta}}{p_{\hat{\boldsymbol{\theta}}\left(\boldsymbol{x}^{n}\right)}\left(\boldsymbol{x}^{n}\right)}=\int_{\Theta}\frac{p_{\boldsymbol{\theta}}\left(\boldsymbol{x}^{n}\right)}{p_{\hat{\boldsymbol{\theta}}\left(\boldsymbol{x}^{n}\right)}\left(\boldsymbol{x}^{n}\right)}\pi\left(\boldsymbol{\theta}\right)\,\mathrm{d}\boldsymbol{\theta}\,,
\]
which proves the lemma.
\end{IEEEproof}
Note that we have not really used that the parametrized family is
an exponential family, so a similar result holds as long as the parametrization
is sufficiently smooth. If the parametrization is sufficiently smooth
one can also prove that the prior must be proportional to Jeffrey's
prior. We conjecture that if conditional MDL is a Bayesian prediction
for some smoothly parametrized family where the parameter space is
finitely dimensional, then the family must be exponential. Recall
that the saddle point approximation \cite{Daniels1954} for the exponential
family is
\[
\exp\left(-nD_{A}\left(\boldsymbol{\theta},\hat{\boldsymbol{\theta}}\left(\boldsymbol{x}^{n}\right)\right)\right)\frac{\left|Cov\left(\boldsymbol{\mu}_{\boldsymbol{\theta}}\right)\right|^{\nicefrac{1}{2}}}{\tau^{\nicefrac{d}{2}}}\,,
\]
where $\tau$ is short for $2\pi.$ 
\begin{thm}
\label{thm:Main}Assume that $\left(P_{\boldsymbol{\theta}}\right)_{\boldsymbol{\theta}\in\Theta}$
is a natural exponential family. Then the following conditions are
equivalent:

$\bullet$ CNML is a Bayesian prediction strategy. 

$\bullet$ Jeffreys' posterior distributions are elements of the conjugated
exponential family.

$\bullet$ The renormalized saddle-point approximation is exact for
the conjugated exponential family.
\end{thm}
\begin{IEEEproof}
According to expression (\ref{lem:key}) we may define a constant
$C_{n}$ by 
\[
C_{n}=\int_{\Theta}\frac{p_{\boldsymbol{\theta}}\left(\boldsymbol{x}^{n}\right)}{p_{\hat{\boldsymbol{\theta}}\left(\boldsymbol{x}^{n}\right)}\left(\boldsymbol{x}^{n}\right)}\pi\left(\boldsymbol{\theta}\right)\,\mathrm{d}\boldsymbol{\theta}.
\]
Then 
\begin{equation}
\frac{p_{\boldsymbol{\theta}}\left(\boldsymbol{x}^{n}\right)}{p_{\hat{\boldsymbol{\theta}}\left(\boldsymbol{x}^{n}\right)}\left(\boldsymbol{x}^{n}\right)}\cdot\frac{\pi\left(\boldsymbol{\theta}\right)}{C_{n}}\label{eq:conj}
\end{equation}
is a probability density function for $\boldsymbol{\theta}$. We will
demonstrate that the family of probability measures (\ref{eq:conj})
parametrized by $\boldsymbol{x}^{n}$ is the conjugated exponential
family with $\boldsymbol{\theta}$ as sufficient statistic. We have
\begin{multline*}
\frac{p_{\boldsymbol{\theta}}\left(\boldsymbol{x}^{n}\right)}{p_{\hat{\boldsymbol{\theta}}\left(\boldsymbol{x}^{n}\right)}\left(\boldsymbol{x}^{n}\right)}\cdot\frac{\pi\left(\boldsymbol{\theta}\right)}{C_{n}}\\
=\frac{\exp\left(n\left(\boldsymbol{\theta}\cdot\bar{\boldsymbol{x}}-A\left(\boldsymbol{\theta}\right)\right)\right)}{\exp\left(n\left(\hat{\theta}\left(\boldsymbol{x}^{n}\right)\cdot\bar{\boldsymbol{x}}-A\left(\hat{\theta}\left(\boldsymbol{x}^{n}\right)\right)\right)\right)}\cdot\frac{\pi\left(\boldsymbol{\theta}\right)}{C_{n}}\\
=\exp\left(n\left(\boldsymbol{\theta}\cdot\bar{\boldsymbol{x}}-A^{*}\left(\bar{\boldsymbol{x}}\right)\right)\right)\cdot\frac{\pi\left(\boldsymbol{\theta}\right)}{\exp\left(nA\left(\boldsymbol{\theta}\right)\right)C_{n}}.
\end{multline*}
According to the robustness property (\ref{eq:robust}) the density
can be rewritten as
\[
\exp\left(-nD_{A}\left(\boldsymbol{\theta},\hat{\boldsymbol{\theta}}\left(\boldsymbol{x}^{n}\right)\right)\right)\cdot\frac{\pi\left(\boldsymbol{\theta}\right)}{C_{n}}.
\]
Since this should hold for $n$ tending to infinity the saddle point
approximation implies that $\pi\left(\boldsymbol{\theta}\right)$
is proportional to $\left|Cov\left(\boldsymbol{\mu}_{\boldsymbol{\theta}}\right)\right|^{\nicefrac{1}{2}}$.
Therefore the density in the exponential family is proportional to
the saddle point approximation.
\end{IEEEproof}
\begin{cor}
If any of the equivalent conditions of Theorem \ref{thm:Main} are
fulfilled the exponential family is steep and the parameter space
is maximal.
\end{cor}
The goal is now to identify exponential families where Jeffreys' posterior
distributions form exponential families with exact renormalized saddle
point approximations. In \cite{Blesild1985} it was proved that under
certain regularity conditions the renormalized saddle point approximation
is exact for \emph{reproductive exponential families}. The reproductive
exponential families were defined and described in detail in \cite{Barndorff-Nielsen1983}
where it was proved in 1 dimension the following families were reproductive:
the Gaussian location families, the Gamma exponential families and
the Inverse Gaussian families. The idea of reproductive exponential
families can be used to construct reproductive exponential families
in higher dimension by combining reproductive exponential families
in lower dimensions. Five non-trivial examples of 2-dimensional (strongly)
reproducible exponential families obtained by combining reproductive
1 dimensional families were listed in \cite{Barndorff-Nielsen1983}.
For each reproductive exponential family the conjugate exponential
family (if it exists) will satisfy the conditions of Theorem \ref{thm:Main}.
We will illustrate how this works for 1-dimensional reproductive exponential
families.

The only 1-dimensional natural exponential families where the renormalized
saddle point approximation is exact, are the three reproductive exponential
families mentioned above \cite{Daniels1980}, and it can be proved
by solving ordinary differential equations \cite{Blesild1985}. A
complete classification of exponential families with exact renormalized
saddle point approximation in dimension 2 or higher would require
solving some complicated partial differential equations. Therefore
a complete catalog of families for which the equivalent conditions
of Theorem \ref{thm:Main} are fulfilled, seems inaccessable.

For the 1-dimensional reproductive exponential families the functions
$A^{*}$ is exactly the ones used in \cite{Barndorff-Nielsen1983}
to prove that the exponential family is reproductive. Exploration
of this fact in higher dimensions will be covered in a future paper. 

\section{The Gamma family}

A Gamma distribution can be parametrized by the shape parameter $\alpha$
and the rate parameter $\beta.$ With these parameters the Gamma distribution
$\Gamma\left(\alpha,\beta\right)$ has density 
\[
\frac{\beta^{\alpha}x^{\alpha-1}}{\Gamma\left(\alpha\right)}\exp\left(-\beta x\right)=\frac{x^{\alpha-1}}{\Gamma\left(\alpha\right)}\exp\left(-\beta x+\alpha\ln\left(\beta\right)\right)
\]
for $x>0.$ For a fixed value of $\alpha$ this is a natural exponential
family with natural parameter $\theta=-\beta<0$. Therefore $A\left(\theta\right)=-\alpha\ln\left(-\theta\right).$
The mean value is $\mu=-\alpha/\theta$ so that $\theta=-\alpha/\mu.$
The variance is $Var=\alpha\theta^{-2}$ , so that the variance function
is $V\left(\mu\right)=\frac{\mu^{2}}{\alpha}.$ In terms of the parameter
$\beta$ the mean value is $\mu=\alpha/\beta$ and the variance is
$Var=\alpha\cdot\beta^{-2}$ . Jeffreys' prior has density proportional
to $\frac{\alpha^{\nicefrac{1}{2}}}{\beta},$ which cannot be normalized. 

The Bregman divergence is
\begin{multline*}
D_{A}\left(\theta_{1},\theta_{2}\right)\\
=\alpha\ln\left(-\frac{1}{\theta_{1}}\right)-\left(\alpha\ln\left(-\frac{1}{\theta_{2}}\right)+\left(\theta_{1}-\theta_{2}\right)\cdot\frac{-\alpha}{\theta_{2}}\right)\\
=\alpha\left(\frac{\theta_{1}}{\theta_{2}}-1-\ln\left(\frac{\theta_{1}}{\theta_{2}}\right)\right).
\end{multline*}
For $\alpha=1$ this Bregman divergence is called the \emph{Itakura-Saito
divergence}. 

The convex conjugate of $A$ is 
\begin{multline*}
A^{*}\left(x\right)=\sup_{\theta}\left\{ x\cdot\theta-A\left(\theta\right)\right\} =x\cdot\left(-\frac{\alpha}{x}\right)-A\left(-\frac{\alpha}{x}\right)\\
=-\alpha+\alpha\ln\left(\frac{\alpha}{x}\right)=-\alpha+\alpha\ln\left(\alpha\right)-\alpha\ln\left(x\right).
\end{multline*}
We see that the conjugated exponential family of $\beta=-\theta$
is again a Gamma exponential family with shape parameter $\alpha$,
i.e. the Gamma exponential family is \emph{self-conjugated}. If $x$
is observed the posterior distribution of $\beta$ has rate parameter
$x$. If a sequence of length $m$ has been observed then the posterior
distribution is a Gamma distribution with shape parameters $m\alpha$
and rate parameter $m\bar{x}.$ 

Since the density of a Gamma distribution equals the re-normalized
saddle point approximation we have that the conditions in Theorem
\ref{thm:Main} are fulfilled and the CNML predictor equals Bayesian
prediction based on Jeffreys' prior. This also holds for exponential
families like the inverse Gamma family, the Pareto family, the Nakagima
family, and the Weibull family where the sufficient statistic is a
smooth 1-to-1 function of the sufficient statistic in a Gamma family.

We will now look at the consequences of self-conjugation for calculations
of one-sided credible intervals and one-sided confidence intervals.

Let $G$ denote the distribution function of $\Gamma\left(m\alpha,m\bar{x}\right)$,
i.e. the posterior distribution of $\beta$ if the average is observed
to be $\bar{x}$. Then $\left[0,G^{-1}\left(1-\tilde{\alpha}\right)\right]$
is a $1-\tilde{\alpha}$ \emph{credible interval} for $\beta.$ We
can write 
\begin{align*}
G^{-1}\left(1-\tilde{\alpha}\right) & =\frac{F^{-1}\left(1-\tilde{\alpha}\right)}{\bar{x}}
\end{align*}
where $F$ is the distribution function of $\Gamma\left(m\alpha,m\right).$
If $X_{i}\sim\Gamma\left(\alpha,\beta\right)$ then $\sum_{i=1}^{m}X_{i}\sim\Gamma\left(m\alpha,\beta\right)$
and $\frac{1}{m}\sum_{i=1}^{m}X_{i}\sim\Gamma\left(m\alpha,m\beta\right)$
so that $\text{\ensuremath{\beta\bar{X}\sim\Gamma\left(m\alpha,m\right)}.}$
Therefore
\begin{align*}
P\left(\beta\in\left[0,\frac{F^{-1}\left(1-\tilde{\alpha}\right)}{\bar{X}}\right]\right) & =P\left(\bar{X}\in\left[0,\frac{F^{-1}\left(1-\tilde{\alpha}\right)}{\beta}\right]\right)\\
 & =1-\tilde{\alpha}
\end{align*}
 so that the $1-\tilde{\alpha}$ credible interval $\left[0,\frac{F^{-1}\left(1-\tilde{\alpha}\right)}{\bar{x}}\right]$
is also a $1-\tilde{\alpha}$ \emph{confidence interval} for $\beta$
as defined in the frequentist approach to statistics. 

\section{The Gaussian location family}

If the parameter space equals $\mathbb{R}^{d}$ the notion of self-conjugation
becomes very simple. The proof of the following lemma is an easy exercise.
\begin{lem}
Let $B:\mathbb{R}^{d}\to\mathbb{R}^{d}$ denote a linear invertible
self-adjoint mapping. If $G$ is a convex function and $F=G\circ B$
then $F^{*}=G^{*}\circ B^{-1}$.
\end{lem}
The Gaussian location model has density 
\[
\frac{\exp\left(-\frac{1}{2}\left(\boldsymbol{x}-\boldsymbol{\mu}\right)\cdot B^{-1}\left(\boldsymbol{x}-\boldsymbol{\mu}\right)\right)}{\tau^{\nicefrac{d}{2}}\cdot\left|B\right|^{^{\nicefrac{1}{2}}}}
\]
where $\boldsymbol{\mu}$ is the mean and $B$ denotes the co-variance
matrix.
\begin{thm}
If an exponential family has a cumulant generating function $A:R^{d}\to R$
that satisfies $A^{*}=A\circ B$ for some positiv definite linear
function $B:\mathbb{R}^{d}\to\mathbb{R}^{d}$ then the exponential
family is a Gaussian location model where $B$ can be identified with
the co-variance matrix.
\end{thm}
\begin{IEEEproof}
Define $F=A\circ B^{\nicefrac{1}{2}}.$ Then 
\begin{multline*}
F^{*}=A^{*}\circ\left(B^{\nicefrac{1}{2}}\right)^{-1}=A\circ B\circ B^{-\nicefrac{1}{2}}=A\circ B^{\nicefrac{1}{2}}=F\,.
\end{multline*}
Since $F$ is self-conjugated and defined on $\mathbb{R}^{d}$ we
can apply \cite[Prop. 29a]{Moreau1965} to get $F\left(\boldsymbol{x}\right)=\frac{1}{2}\left\Vert \boldsymbol{x}\right\Vert ^{2}.$
Therefore 
\begin{multline*}
A\left(\boldsymbol{x}\right)=F\left(B^{-\nicefrac{1}{2}}\left(\boldsymbol{x}\right)\right)=\frac{1}{2}B^{-\nicefrac{1}{2}}\left(\boldsymbol{x}\right)\cdot B^{-\nicefrac{1}{2}}\left(\boldsymbol{x}\right)\\
=\frac{1}{2}\boldsymbol{x}\cdot B^{-1}\left(\boldsymbol{x}\right).
\end{multline*}
It is easy to prove that the Gaussian location model also has cumulant
generating function $\frac{1}{2}\boldsymbol{x}\cdot B^{-1}\left(\boldsymbol{x}\right).$
\end{IEEEproof}
Since the saddle point approximation is exact for the Gaussian location
family the conditions of Threorem \ref{thm:Main} are fulfilled. For
the Gaussian location family the Bregman divergence is symmetric in
its arguments and inference reduces to the principle of least squares. 

In Bayesian statistics a $1-\tilde{\alpha}$ \emph{credible region}
for the mean value parameter can be calculated as a divergence ball
\begin{equation}
\left\{ \theta\in R^{d}\mid D_{A}\left(\boldsymbol{\theta},\hat{\boldsymbol{\theta}}\left(\boldsymbol{x}\right)\right)\leq r\right\} \label{eq:ball}
\end{equation}
where the radius $r$ is chosen so that the ball has probability $1-\tilde{\alpha}$.
Using that the exponential family is self-conjugated we see that the
ball (\ref{eq:ball}) is also a $1-\tilde{\alpha}$ \emph{confidence
region} as defined in frequentist statistics.

\section{The Poisson-exponential family}

The saddle point approximation is exact for the inverse Gaussian family
with density 
\[
\left(\frac{\kappa}{\tau\beta^{3}}\right)^{\nicefrac{1}{2}}\exp\left(-\kappa\frac{\left(\beta-\beta_{0}\right)^{2}}{2\beta_{0}^{2}\beta}\right),
\]
where $\beta$ is the sufficient statistic and $\beta_{0}$ denotes
the mean value of the distribution and $\kappa$ denotes the \emph{shape
parameter}. We are going to identify the conjugated exponential family.
First we rewrite
\begin{multline*}
\left(\frac{\kappa}{\tau\beta^{3}}\right)^{\nicefrac{1}{2}}\exp\left(-\kappa\frac{\left(\beta-\beta_{0}\right)^{2}}{2\beta_{0}^{2}\beta}\right)\\
=\left(\frac{\kappa}{\tau\beta^{3}}\right)^{\nicefrac{1}{2}}\exp\left(-\frac{\kappa}{2\beta}\right)\exp\left(\frac{-\kappa}{2\beta_{0}^{2}}\cdot\beta+\frac{\kappa}{\beta_{0}}\right).
\end{multline*}
The natural parameter is $\theta=\frac{-\kappa}{2\beta_{0}^{2}}$
and the cumulant generating function is $A\left(\theta\right)=\left(-2\kappa\theta\right)^{\nicefrac{1}{2}}.$ 

The convex conjugate is 
\begin{align*}
A^{^{*}}\left(\beta\right) & =\sup\left\{ \beta\cdot\theta-A\left(\theta\right)\right\} \\
 & =\beta\cdot\frac{-\kappa}{2\beta^{2}}-\left(-2\kappa\cdot\frac{-\kappa}{2\beta^{2}}\right)^{\nicefrac{1}{2}}=\frac{\kappa}{2\beta}.
\end{align*}
One can identify an exponential family with this function as cumulant
generating function by taking the inverse Laplace transform, but it
is more instructive to identify it by calculating the variance function.
We have 
\begin{align*}
\left(A^{*}\right)'\left(\beta\right) & =-\frac{\kappa}{2\beta^{2}}\,\textrm{and}\,\left(A^{*}\right)''\left(\beta\right)=\frac{\kappa}{\beta^{3}}.
\end{align*}
Thus $\hat{\theta}\left(\beta\right)=-\frac{\kappa}{2}\beta^{-2}$
so that $\beta\left(\theta\right)=\left(-\frac{\kappa}{2\theta}\right)^{\nicefrac{1}{2}}$
and $V\left(\theta\right)=\kappa\left(\beta\left(\theta\right)\right)^{-3}=2^{\nicefrac{3}{2}}\kappa^{-\nicefrac{1}{2}}\left(-\theta\right)^{\nicefrac{3}{2}}=\phi\cdot\left(-\theta\right)^{\nicefrac{3}{2}}$
where $\phi=2^{\nicefrac{3}{2}}\kappa^{-\nicefrac{1}{2}}$. Since
the variance function is a power function of order $\nicefrac{3}{2}$
one says that the corresponding exponential family is a \emph{Tweedie
family} of order $p=\nicefrac{3}{2}$~. Jeffreys' prior for this
family is proportional to 
\[
\left(\left(A^{*}\right)''\left(\beta\right)\right)^{\nicefrac{1}{2}}=\kappa^{\nicefrac{1}{2}}\cdot\beta^{-\nicefrac{3}{2}},
\]
which cannot be normalized. Credible intervals and confidence intervals
can be calculated using $\mathtt{tweedie}$ and the $\mathtt{statmod}$
package in the R program, but the $1-\tilde{\alpha}$ credible intervals
do not coincide with the $1-\tilde{\alpha}$ confidence intervals
reflecting that the Poisson-exponential family is not self-conjugated.

One cannot calculate the density of elements of the Tweedie family
of order $p=\nicefrac{3}{2}$ exactly, but they can be obtained by
the following construction. Let $N$ denote a random variable with
a Poisson distribution $Po\left(\lambda\right)$. Let $X_{1},X_{2},\dots$
denote a sequence of iid random variables each exponentially distributed
$Exp\left(\beta\right)$. Then we may define 
\[
Y=\sum_{n=1}^{N}X_{n}\,.
\]
Then the distribution of $Y$ is a compound Poisson distribution.
Distributions where $X_{i}$ are Gamma distributions were called Poisson-gamma
distributions in \cite{Smyth1996}, so we will call the distribution
of $Y$ a \emph{Poisson-exponential distribution} when $X_{i}$ are
exponential. The density of $\sum_{n=1}^{\alpha}X_{n}$ is 
\[
\frac{\tilde{\beta}^{\alpha}x^{\alpha-1}\exp\left(-\tilde{\beta}x\right)}{\Gamma\left(\alpha\right)}.
\]
Therefore the Poisson-exponential distribution has a point mass in
0 of weight $\exp\left(-\lambda\right)$ and it has density
\[
\sum_{\alpha=0}^{\infty}\frac{\lambda^{\alpha}\exp\left(-\lambda\right)}{\alpha!}\cdot\frac{\beta^{\alpha}x^{\alpha-1}\exp\left(-\beta x\right)}{\Gamma\left(\alpha\right)}
\]
for $x>0$. We introduce $\kappa=\frac{\tilde{\beta}\cdot\lambda}{2}$
so that the density can be written as
\[
\sum_{\alpha=0}^{\infty}\frac{\left(\frac{\kappa}{2}\right)^{\alpha}x^{\alpha-1}}{\alpha!\Gamma\left(\alpha\right)}\cdot\exp\left(-\beta\cdot x-\frac{\kappa}{2\beta}\right).
\]
This is a natural exponential family with with natural parameter $-\beta$
and cumulant generating function $\kappa/\left(2\beta\right)$. Except
for a change of sign it is the conjugated exponential family of the
inverse Gaussian family.

Since the saddle point approximation is exact for the inverse Gaussian
family, prediction for the Poisson-exponential family based on CNML
equals prediction based on Jeffreys prior, and Jeffreys posterior
equals an inverse Gaussian distribution. 

The Poisson-exponential families have been used to model the accumulated
amount of rain in rainfalls, where the amount of rain in each rainfall
is modeled by an exponential distribution and the number of rainfalls
is modeled by a Poisson distribution \cite{Thompson1984,Revfeim1984}.
This application dates back to Cornish and Fisher. Reference to other
applications as well as a derivation of the basic properties of Poisson-gamma
distributions can be found in \cite{Withers2011}. Note that the Poisson-exponential
family is a Tweedie family of order $p=\nicefrac{3}{2}$ and that
some of the literature on applications of the Poisson-exponential
family treat the order $p$ as a free parameter that should be estimated
in order to give a good fit with data. According to our results the
value $p=\nicefrac{3}{2}$ is special with respect to statistical
inference, so that $p$ cannot be considered as a free parameter if
we want to have the properties developed here.

\section*{Acknowledgement}

I would like to thank Wojciech Kot\l owski for useful comments to
this paper.

{\small{}

}{\small\par}
\end{document}